\documentclass[12pt]{amsart}
\usepackage{latexsym}
\usepackage{amssymb} 
\usepackage{mathrsfs}
\usepackage{amsmath}
\usepackage{latexsym}
\usepackage{delarray}
\usepackage{amssymb,amsmath,amsfonts,amsthm,
mathrsfs}

\usepackage[colorlinks]{hyperref}

\renewcommand\eqref[1]{(\ref{#1})} 

\newtheorem{theorem}{Theorem}[section]

\newtheorem{proposition}[theorem]{Proposition}

\theoremstyle{definition}
\newtheorem{remark}[theorem]{Remark}

\newcommand{\Cinfc}{\ensuremath{\mathcal{C}^\infty_{\text{c}}}}

\newcommand{\E}{\ensuremath{{\mathcal E}}}

\newcommand{\mb}[1]{\ensuremath{\mathbb{#1}}}

\newcommand{\R}{\mb{R}}

\newcommand{\lara}[1]{\langle #1 \rangle}

\newfont{\bl}{msbm10 scaled \magstep2}

\newcommand{\beq}{\begin{equation}}
\newcommand{\eeq}{\end{equation}}

\newcommand{\eps}{\varepsilon}

\newcommand{\esp}{\mathrm{e}}
\newcommand\Rn{{\mathbb R}^n}
\newcommand{\irm}{{\rm i}}

\title[Hyperbolic systems with H\"older characteristics]{
On hyperbolic systems with time dependent H\"older characteristics}

\author[Claudia Garetto]{Claudia Garetto}
\address{
  Claudia Garetto:
  \endgraf
  Department of Mathematical Sciences
  \endgraf
  Loughborough University
  \endgraf
  Loughborough, Leicestershire, LE11 3TU
  \endgraf
  United Kingdom
  \endgraf
  {\it E-mail address} {\rm c.garetto@lboro.ac.uk}
  }
\author[Michael Ruzhansky]{Michael Ruzhansky}
\address{
 Michael Ruzhansky:
  \endgraf
  Department of Mathematics
  \endgraf
  Imperial College London
  \endgraf
  180 Queen's Gate, London SW7 2AZ
  \endgraf
  United Kingdom
  \endgraf
  {\it E-mail address} {\rm m.ruzhansky@imperial.ac.uk}
  }

\thanks{The first author was supported by the
EPSRC First grant EP/L026422/1.
The second author was supported by the
EPSRC grant EP/E062873/1.
No new data was collected or generated during the course of the research.
}
\date{}


\begin{document}

\maketitle

\begin{abstract}
In this paper we study the well-posedness of weakly hyperbolic systems with time dependent coefficients. We assume that the eigenvalues are low regular, in the sense that they are H\"older with respect to $t$. In the past these kind of systems have been investigated by Yuzawa \cite{Yu:05} and Kajitani \cite{KY:06} by employing semigroup techniques (Tanabe-Sobolevski method). Here, under a certain uniform property of the eigenvalues, we improve the Gevrey well-posedness result of \cite{Yu:05} and we obtain well-posedness in spaces of ultradistributions as well. Our main idea is a reduction of the system to block Sylvester form and then the formulation of suitable energy estimates inspired by the treatment of scalar equations in \cite{GR:11}.   \end{abstract}

\section{Introduction}
We want to study the Cauchy problem for first order hyperbolic systems of the type
\beq
\label{CP_syst_fin}
\begin{split}
D_t u-A(t,D_x)u-B(t)u&=0,\quad x\in\R^n,\, t\in[0,T],\\
u|_{t=0}&=g_0,
\end{split}
\eeq
where $A$ and $B$ are $m\times m$ matrices of first order and zero order differential operators, respectively, with $t$-dependent coefficients, $u$ and $g_0$ are column vectors with $m$ entries. We work under the assumptions that the system matrix is of size $m\times m$ with real eigenvalues and that the coefficients are of class $C^{m-1}$ with respect to $t$. It follows that at the points of highest multiplicity the eigenvalues are of H\"older class $(m-1)/m$. We will therefore assume that the matrix $A(t,\xi)$ has $m$ real eigenvalues $\lambda_j(t,\xi)$ of H\"older class $C^\alpha$, $0<\alpha\le1$ with respect to $t$.  Note that it is not restrictive to assume that the eigenvalues $\lambda_j$, $j=1,\dots,m$, are ordered because
we can always reorder them to satisfy this (ordering) assumption, and the H\"older continuity is preserved by such
reordering. If $\alpha=1$, it is sufficient to assume that $\lambda_j$, $j=1,\dots,m$, is Lipschitz.

In analogy with scalar equations in \cite{GR:11} and \cite{ColKi:02}
we work under the hypothesis of the following \emph{uniform property}: there exists a constant $c>0$ such that 
\beq
\label{hyp_coincide}
|\lambda_i(t,\xi)-\lambda_j(t,\xi)|\le c|\lambda_k(t,\xi)-\lambda_{k-1}(t,\xi)|,
\eeq
for all $1\le i,j,k\le m$, $t\in[0,T]$ and $\xi\in\R^n$. 






Assumptions of H\"older regularity of this type and the uniform condition \eqref{hyp_coincide}
are rather natural, see Colombini and Kinoshita \cite{ColKi:02} and
the authors' paper \cite{GR:11} for a discussion and examples.
In particular, Colombini and Kinoshita \cite{ColKi:02} treated the scalar version of the Cauchy problem
\eqref{CP_syst_fin} with $n=1$, and the authors extended it to the multidimensional case $n\geq 1$
in \cite{GR:11}, also improving some Gevrey indices.

The research of this paper continues investigations of properties of solutions to Cauchy problems 
for hyperbolic equations with multiplicities. The case of time-dependent coefficients already presents
a number of challenging problems, most importantly in view of the fact that already the scalar wave equation
\beq\label{WE}
\partial_t^2 v-a(t)\Delta v=0,\; v(0)=v_0,\; \partial_t v(0)=v_1,
\eeq
in dimension $n=1$ may not be well-posed even for smooth data $v_0,v_1\in C^\infty$.
More precisely, if $a\in C^\alpha$ is H\"older with $0<\alpha<1$, even in the strictly hyperbolic case 
$a>0$ the Cauchy problem \eqref{WE} may have non-unique solutions
(see Colombini, Jannelli and Spagnolo \cite{Colombini-Jannelli-Spagnolo:Annals-low-reg});
and in the weakly hyperbolic case $a\geq 0$ even if $a$ is smooth $a\in C^\infty$, 
the Cauchy problem \eqref{WE} may have no distributional solutions
(see Colombini and Spagnolo \cite{Colombini-Spagnolo:Acta-ex-weakly-hyp}).
However, the Cauchy problem \eqref{WE} is well-posed in suitable Gevrey classes,
see Colombini, de Giorgi and Spagnolo \cite{Colombini-deGiordi-Spagnolo-Pisa-1979}.
At the moment, scalar higher order equations with time-dependent coefficients are relatively 
understood, see e.g. \cite{ColKi:02,KS} and their respective extensions in
\cite{GR:11,GR:12}. Further extreme cases: analytic coefficients and distributional coefficients have
been also investigated, see e.g. authors' papers \cite{Garetto-Ruzhansky:JMAA-2014,Garetto-Ruzhansky:ARMA, G:15},
respectively, and references therein. Hyperbolic systems of the form \eqref{CP_syst_fin} have been also investigated, see e.g.
Yuzawa \cite{Yu:05}, Kajitani and Yuzawa \cite{KY:06} and Garetto \cite{G:15}.

\medskip
The main new idea behind this paper enabling us to obtain an improvement in the well-posedness results
for the system in \eqref{CP_syst_fin}
is the transformation of the system \eqref{CP_syst_fin} to a larger system which, however, enjoys the property of
being in block Sylvester form. Such a transformation, which can be performed under the assumption that the system coefficients are of class $C^{m-1}$ with respect to $t$, is carried out following the method of D'Ancona and Spagnolo
\cite{DS}, leading to the Cauchy problem of the form
\beq
\label{CP_syst_Syl0}
\begin{split}
D_tU-\mathcal{A}(t,D_x)U-\mathcal{L}(t,D_x)U&=0,\\
U_{t=0} &= \{D_t^{j-1}\lara{D_x}^{m-j}g_0\}_{j=1,2,\dots,m}.
\end{split}
\eeq
This is a Cauchy problem for the first order hyperbolic system of the size $m^2\times m^2$ of
pseudo-differential equations. Despite the increase of the size of the system from $m\times m$ to $m^2\times m^2$ and
the change from a differential system to a pseudo-differential one, the system
\eqref{CP_syst_Syl0} has a crucial advantage of being in a block Sylvester form, see
\eqref{CP_syst_Syl} for a precise formulation. This allows us to implement the ideas developed in
\cite{GR:11} for scalar equations, where the reduction of a scalar equation to a Sylvester form system was performed.
 
To summarise our result here we first note that combining the results in \cite{Yu:05, KY:06} we already know that the Cauchy problem \eqref{CP_syst_fin} is well-posed in the Gevrey class $\gamma^s$, with 
\beq
\label{int_yu}
1\le s<1+\frac{\alpha}{m}.
\eeq
Arguing by the Fourier characterisation of Gevrey-Beurling ultradistributions one can easily extend the Gevrey well-posedness above to spaces of ultradistributions. It is our aim in this paper to show that the interval of Gevrey well-posedness in \eqref{int_yu} can be enlarged under the uniform property \eqref{hyp_coincide} of the eigenvalues. Since by the results of Kajitani and Yuzawa at least an ultradistributional solution exists for Gevrey initial data with $s\ge 1+\frac{\alpha}{m}$ we will prove, for suitable values of $s$, that this solution is indeed Gevrey, because it solves the reduced Cauchy problem \eqref{CP_syst_Syl0}. In this sense, the well-posedness of \eqref{CP_syst_fin} can be determined by studying the well-posedness of the reduced Cauchy problem  \eqref{CP_syst_Syl0}. More precisely, by standard arguments it is sufficient to find an a-priori estimate on the Fourier transform with respect to $x$ of the solution $U$ of \eqref{CP_syst_Syl0}.

We assume that the Gevrey classes $\gamma^s(\Rn)$ are well-known: these are spaces of all
 $f\in C^\infty(\R^n)$ such that for every compact set $K\subset\R^n$ there
exists a constant $C>0$ such that for all $\beta\in\mathbb N_0^n$ we have the estimate
\beq\label{Gevrey1}
\sup_{x\in K}|\partial^\beta f(x)|\le C^{|\beta|+1}(\beta!)^s.
\eeq
For $s=1$, we obtain the class of analytic functions.
We refer to 
\cite{GR:11} for a detailed discussion and Fourier characterisations of Gevrey spaces of different types. 
Since we are dealing with vectors in
this paper, we will write $\gamma^s(\Rn)^m$ for $m$-vectors consisting of functions in $\gamma^s(\Rn)$.
This is our main result:

\begin{theorem}
\label{thm1}
Assume that coefficients of the $m\times m$ matrices $A$ and $B$ are of class $C^{m-1}$ and that the matrix $A(t,\xi)$ has $m$ real eigenvalues $\lambda_j(t,\xi)$ of H\"older class $C^\alpha$, $0<\alpha\le1$ with respect to $t$, that satisfy
\eqref{hyp_coincide}. Let $T>0$ and $g_0\in\gamma^s(\Rn)^m$.
Then,  the Cauchy problem \eqref{CP_syst_fin} has a unique solution
$u\in C^1([0,T],\gamma^s(\Rn)^m)$ provided that
\beq
\label{s1}
1\leq s< 1+\min\biggl\{\alpha,\frac{1}{m-1}\biggr\}.
\eeq

\end{theorem}

For the proof we can assume that $s>1$ since the case $s=1$ is essentially known, see
\cite{Jannelli:analytic-CPDE-1984} and \cite{Kajitani:analytic-CPDE-1986}.

Also, we note that the proof also covers the case $\alpha=1$, in which case it is enough to
assume that the eigenvalues are Lipschitz.

We note that the result of Theorem \ref{thm1} is an improvement of known results in terms of the Gevrey order.
For example, this is an improvement of Yuzawa's  and Kajitani's order \eqref{int_yu} from \cite{KY:06, Yu:05}. See Remark \ref{rem_Yu} for more details.

\medskip
The energy estimates obtained in the proof of Theorem \ref{thm1} allow one to also obtain the 
ultradistributional well-posedness results. First we note that the Gevrey spaces
$\gamma^s(\R^n)$ considered in \eqref{Gevrey1} are of Gevrey-Roumeau type. At the same time, we denote by
$\gamma^{(s)}(\R^n)$ the Gevrey spaces of Gevrey-Beurling type,
i.e. the space of all $f\in C^\infty(\Rn)$ such that
for every compact set $K\subset\R^n$
and for every constant $A>0$ there exists a constant $C_{A,K}>0$ such that for all $\beta\in\mathbb N_0^n$ we have the estimate
\[
\sup_{x\in K}|\partial^\beta f(x)|\le C_{A,K} A^{|\beta|}(\beta!)^s.
\]
For $1<s<\infty$, we denote by
$\mathcal D_{(s)}^\prime(\R^n):=(\gamma^{(s)}_c(\R^n))'$ the topological dual of compactly supported functions in $\gamma^{(s)}(\R^n)$ and by $\E'_{(s)}(\R^n)$ the topological dual of $\gamma^{(s)}(\R^n)$. 
Consequently, arguing similarly to \cite{GR:11}, the proof of Theorem \ref{thm1} yields the
following ultradistributional well-posedness:

\begin{theorem}
\label{thm2}
Assume that coefficients of the $m\times m$ matrices $A$ and $B$ are of class $C^{m-1}$ and that the matrix $A(t,\xi)$ has $m$ real eigenvalues $\lambda_j(t,\xi)$ of H\"older class $C^\alpha$, $0<\alpha\le1$ with respect to $t$, that satisfy
\eqref{hyp_coincide}. Let $T>0$ and $g_0\in (\E_{(s)}^\prime(\R^n))^m$.
Then, the Cauchy problem \eqref{CP_syst_fin} has a unique solution
$u\in C^1([0,T],(\mathcal D_{(s)}^\prime(\R^n))^m)$ provided that
$$
\label{s1d}
1<s\leq 1+\min\biggl\{\alpha,\frac{1}{m-1}\biggr\}.
$$

\end{theorem}

\section{Proof of Theorem \ref{thm1}}

The first step in our new approach to the Cauchy problem \eqref{CP_syst_fin} is to rewrite the system in a special form, i.e., in block Sylvester form. This is possible thanks to the reduction given by d'Ancona and Spagnolo in \cite{DS}, which is summarised in the following subsection.

\subsection{Reduction to block Sylvester form}

We begin by considering the cofactor matrix $L(t,\tau,\xi)$ of $(\tau I-A(t,\xi))^T$ where $I$ is the $m\times m$ identity matrix. By applying the corresponding operator $L(t,D_t,D_x)$ to \eqref{CP_syst_fin} we transform the system
\[
D_t u-A(t,D_x)u-B(t)u=0
\]
into
\beq
\label{red_1}
\delta(t,D_t,D_x)Iu-C(t,D_t,D_x)u=0,
\eeq
where $\delta(t,\tau,\xi)={\rm det}(\tau I-A(t,\xi))$, $C(t,D_t,D_x)$ is the matrix of  lower order terms (differential operators of order $m-1$). Since the entries of $A$ and $B$ are of class $C^{m-1}$ with respect to $t$ the equation above has continuous $t$-dependent coefficients. Indeed, the coefficients of the equation $D_t u-A(t,D_x)u-B(t)u=0$ are of class $C^{m-1}$ and the operator $L(t,D_t,D_x)$ is of order $m-1$ being defined via the cofactor matrix of a $m\times m$ matrix. Note that $\delta(t,D_t,D_x)$ is the operator
\[
D_t^m+\sum_{h=0}^{m-1}b_{m-h}(t,D_x)D_t^h,
\]
with $b_{m-h}(t,\xi)$ homogeneous polynomial of order $m-h$.

We got in this way a set of scalar equations of order $m$ which can be transformed into a first order system of size $m^2\times m^2$ of pseudodifferential equations, by setting
\[
U=\{D_t^{j-1}\lara{D_x}^{m-j}u\}_{j=1,2,\dots,m},
\]
where $\lara{D_x}$ is the pseudodifferential operator with symbol $\lara{\xi}$. More precisely, the equation \eqref{red_1} is now written as
\[
D_tU-\mathcal{A}(t,D_x)U-\mathcal{L}(t,D_x)U=0,
\]
where $\mathcal{A}$ is a $m^2\times m^2$ matrix made of $m$ identical blocks of the type
\begin{multline*}
\lara{D_x}\cdot\\
\left(\begin{array}{cccccc}
                              0 & 1 & 0 & \cdots & 0 \\
                               0 & 0 & 1 &  \cdots & 0\\
                               \vdots & \vdots & \vdots & \vdots & \vdots\\
                              -b_m(t,D_x)\lara{D_x}^{-m}& -b_{m-1}(t,D_x)\lara{D_x}^{-m+1} & \cdots  & \cdots & -b_1(t,D_x)\lara{D_x}^{-1}\\
                              \end{array}
                           \right),
                           \end{multline*}
and the matrix $\mathcal{L}$  of the lower order terms is made of $m$ blocks of  size $m\times m^2$ of the type
\[
\left(\begin{array}{cccccc}
                              0 & 0 & 0 & \cdots & 0 & 0 \\
                               0 & 0 & 0 &  \cdots &  0 & 0\\
                               \vdots & \vdots & \vdots & \vdots & \vdots\\
                              l_{j,1}(t, D_x)& l_{j,2}(t,D_x) & \cdots  & \cdots & l_{j,m^2-1}(t,D_x) & l_{j,m^2}(t,D_x)
                                                     \end{array}
                           \right), 
\]
with $j=1,\dots,m$.
Note that the entries of the matrices $\mathcal{A}$ and $\mathcal{L}$ are pseudodifferential operators of order $1$ and $0$, respectively. 

Concluding, the Cauchy problem \eqref{CP_syst_fin} has been transformed into
\beq
\label{CP_syst_Syl}
\begin{split}
D_tU-\mathcal{A}(t,D_x)U-\mathcal{L}(t,D_x)U&=0,\\
U_{t=0} &= \{D_t^{j-1}\lara{D_x}^{m-j}g_0\}_{j=1,2,\dots,m}.
\end{split}
\eeq
This is a Cauchy problem of first order pseudodifferential equations with principal part in block Sylvester form. The size of the system is increased from $m\times m$ to $m^2\times m^2$ but the system is still hyperbolic, since the eigenvalues of any block of $\mathcal{A}(t,\xi)$ are the eigenvalues of the matrix $A(t,\xi)$.



\subsection{Energy estimates}
As in \cite{GR:11} we regularise the eigenvalues $\lambda_{j}(t,\xi)$ with respect to $t$ and we separate them by adding some power of a parameter $\eps\to 0$. In detail, assuming that the $\lambda_j$'s are ordered and taking a mollifier $\varphi\in\Cinfc(\R)$, $\varphi\ge 0$ with $\int\varphi(t)\, dt=1$ we set
\[
\lambda_{j,\eps}(t,\xi):=(\lambda_j(\cdot,\xi)\ast\varphi_\eps)(t)+j\eps^\alpha\lara{\xi},\qquad t\in[0,T],\, \xi\in\R^n,
\]
where $\varphi_\eps(t)=\eps^{-1}\varphi(t/\eps)$ and $j=1,\dots,m$. The next proposition collects the main properties of these regularised eigenvalues and has been proven in \cite{GR:11} (see Propositions 18 and 19). 
\begin{proposition}
\label{prop_roots}
Let $\varphi\in C^\infty_{c}(\R)$, $\varphi\ge 0$ with $\int_\R\varphi(x)\, dx=1$.

\medskip
Under the assumptions of Theorem \ref{thm1}, let
\beq
\label{def_lambdaj_1}
\lambda_{j,\eps}(t,\xi):=(\lambda_j(\cdot,\xi)\ast\varphi_\eps)(t)+j\eps^\alpha\lara{\xi},
\eeq
for $j=1,...,m$ and $\varphi_\eps(s)=\eps^{-1}\varphi(s/\eps)$, $\eps>0$.
Then, there exists a constant $c>0$ such that
\begin{itemize}
\item[(i)] $|\partial_t\lambda_{j,\eps}(t,\xi)|\le c\,\eps^{\alpha-1}\lara{\xi}$,
\item[(ii)] $|\lambda_{j,\eps}(t,\xi)-\lambda_j(t,\xi)|\le c\,\eps^{\alpha}\lara{\xi}$,
\item[(iii)] $\lambda_{j,\eps}(t,\xi)-\lambda_{i}(t,\xi)\ge \eps^\alpha\lara{\xi}$ for $j>i$,
\end{itemize}
for all $t,s\in[0,T']$ with $T'<T$ and all $\xi\in\R^n$.\\

\end{proposition}

We can now define the $m^2\times m^2$ block diagonal matrix $H_\eps$ made of $m$ identical blocks of the type
\beq
\label{def_block_H}
\left(
    \begin{array}{ccccc}
      1 & 1 & 1 & \dots & 1\\
      \lambda_{1,\eps}\lara{\xi}^{-1} & \lambda_{2,\eps}\lara{\xi}^{-1} & \lambda_{3,\eps}\lara{\xi}^{-1} & \dots & \lambda_{m,\eps}\lara{\xi}^{-1} \\
      \lambda^2_{1,\eps}\lara{\xi}^{-2} & \lambda^2_{2,\eps}\lara{\xi}^{-2} & \lambda^2_{3,\eps}\lara{\xi}^{-2} & \dots & \lambda^2_{m,\eps}\lara{\xi}^{-2} \\
      \dots & \dots & \dots & \dots & \dots\\
      \lambda^{m-1}_{1,\eps}\lara{\xi}^{-m+1} & \lambda^{m-1}_{2,\eps}\lara{\xi}^{-m+1} & \lambda^{m-1}_{3,\eps}\lara{\xi}^{-m+1} & \dots & \lambda^{m-1}_{m,\eps}\lara{\xi}^{-m+1}\\
    \end{array}
  \right).
\eeq
By separation of the regularised eigenvalues one easily sees that the matrix $H_\eps$ is invertible. Since weakly hyperbolic equations and systems posses the finite speed of propagation property, we know that if the initial data is compactly supported then the solution will be compactly supported in $x$ as well. Hence, instead of dealing with the Cauchy problem \eqref{CP_syst_Syl} directly we can apply the Fourier transform with respect to $x$ to  it and focus on the corresponding Cauchy problem
\beq
\label{CP_syst_Syl_F}
\begin{split}
D_tV-\mathcal{A}(t,\xi)V-\mathcal{L}(t,\xi)V&=0,\\
V_{t=0} &= \{D_t^{j-1}\lara{\xi}^{m-j}\widehat{g_0}\}_{j=1,2,\dots,m}.
\end{split}
\eeq
Note assuming compactly supported initial data in Theorem \ref{thm1} is not restrictive. We look for a solution $V(t,\xi)$ of the type
\beq
\label{VW}
V(t,\xi)=\esp^{-\rho(t)\lara{\xi}^{\frac{1}{s}}}(\det H_\eps)^{-1}H_\eps W,
\eeq
where $\rho\in C^1[0,T]$ will be determined in the sequel. By substitution in \eqref{CP_syst_Syl_F} we obtain
\begin{multline*}
\esp^{-\rho(t)\lara{\xi}^{\frac{1}{s}}}(\det H_\eps)^{-1}H_\eps D_tW+\esp^{-\rho(t)\lara{\xi}^{\frac{1}{s}}}\irm\rho'(t)\lara{\xi}^{\frac{1}{s}}(\det H_\eps)^{-1}H_\eps W+\\
+\irm\esp^{-\rho(t)\lara{\xi}^{\frac{1}{s}}}\frac{\partial_t\det H_\eps}{(\det H_\eps)^2}H_\eps W +\esp^{-\rho(t)\lara{\xi}^{\frac{1}{s}}}(\det H_\eps)^{-1}(D_tH_\eps)W\\
=\esp^{-\rho(t)\lara{\xi}^{\frac{1}{s}}}(\det H_\eps)^{-1}(\mathcal{A}+\mathcal{L})H_\eps W.
\end{multline*}
Multiplying both sides of the previous equation by $\esp^{\rho(t)\lara{\xi}^{\frac{1}{s}}}(\det H_\eps)H_\eps^{-1}$ we get
\[
D_tW+\irm\rho'(t)\lara{\xi}^{\frac{1}{s}}W+\irm\frac{\partial_t\det H_\eps}{\det H_\eps}W + H_\eps^{-1}(D_t H_\eps)W= H_\eps^{-1}(\mathcal{A}+\mathcal{L})H_\eps W.
\]
Thus,
\begin{multline}
\label{energy}
\partial_t |W(t,\xi)|^2=2{\rm Re} (\partial_t W(t,\xi),W(t,\xi))\\
=2\rho'(t)\lara{\xi}^{\frac{1}{s}}|W(t,\xi)|^2+2\frac{\partial_t\det H_\eps}{\det H_\eps}|W(t,\xi)|^2-2
{\rm Re}(H_\eps^{-1}\partial_t H_\eps W,W)\\
-2{\rm Im} (H_\eps^{-1}\mathcal{A}H_\eps W,W)-2{\rm Im} (H_\eps^{-1}\mathcal{L}H_\eps W,W).
\end{multline}
Inspired by the treatment of higher order equations given in \cite{GR:11} we proceed by estimating the terms:
\begin{enumerate}
\item $\frac{\partial_t\det H_\eps}{\det H_\eps}$,
\item $\Vert H_\eps^{-1}\partial_t H_\eps\Vert$,
\item $\Vert H_\eps^{-1}\mathcal{A}H_\eps-(H_\eps^{-1}\mathcal{A}H_\eps)^\ast\Vert$,
\item $\Vert H_\eps^{-1}\mathcal{L}H_\eps-(H_\eps^{-1}\mathcal{L}H_\eps)^\ast\Vert$.
\end{enumerate}

\subsubsection{{\bf Estimate of (i), (ii), (iii) and (iv)}}

We begin by noting that the $m$ identical blocks of the $m^2\times m^2$-matrix $H_\eps$ are exactly given by the matrix $H$ used in the paper 
\cite{GR:11} (formula (3.4)). Hence we can set
\[
H=\left(
    \begin{array}{ccccc}
      1 & 1 & 1 & \dots & 1\\
      \lambda_{1,\eps}\lara{\xi}^{-1} & \lambda_{2,\eps}\lara{\xi}^{-1} & \lambda_{3,\eps}\lara{\xi}^{-1} & \dots & \lambda_{m,\eps}\lara{\xi}^{-1} \\
      \lambda^2_{1,\eps}\lara{\xi}^{-2} & \lambda^2_{2,\eps}\lara{\xi}^{-2} & \lambda^2_{3,\eps}\lara{\xi}^{-2} & \dots & \lambda^2_{m,\eps}\lara{\xi}^{-2} \\
      \dots & \dots & \dots & \dots & \dots\\
      \lambda^{m-1}_{1,\eps}\lara{\xi}^{-m+1} & \lambda^{m-1}_{2,\eps}\lara{\xi}^{-m+1} & \lambda^{m-1}_{3,\eps}\lara{\xi}^{-m+1} & \dots & \lambda^{m-1}_{m,\eps}\lara{\xi}^{-m+1}\\
    \end{array}
  \right).
\]
and observe that
\[
\frac{\partial_t\det H_\eps}{\det H_\eps}=\frac{\partial_t\det H}{\det H}.
\]
By arguing as in (4.3) in \cite{GR:11} we immediately have that 
\beq
\label{est_1_c1}
\biggl|\frac{\partial_t\det H_\eps(t,\xi)}{\det H_\eps(t,\xi)}\biggr|\le c_1\eps^{-1},
\eeq
for all $t\in[0,T]$, $\xi\in\R^n$ and $\eps\in(0,1]$.

Since $H_\eps$ is block diagonal its inverse will be block diagonal as well and precisely given by $m$ identical blocks $H^{-1}$ as defined in Proposition 17(ii) in \cite{GR:11}. It follows that to estimate  $\Vert H_\eps^{-1}\partial_t H_\eps\Vert$ it is enough to estimate the norm of the corresponding block  $H^{-1}\partial_t H$. This has been done in Subsection 4.2 in \cite{GR:11} and leads to
\beq
\label{est_2_c1}
\Vert H_\eps^{-1}\partial_t H_\eps\Vert\le c_2\eps^{-1}.
\eeq
Note that to obtain \eqref{est_2_c1} one uses the uniform property \eqref{hyp_coincide} of the eigenvalues and of the corresponding regularisations.

The same block argument applies to $\Vert H_\eps^{-1}\mathcal{A}H_\eps-(H_\eps^{-1}\mathcal{A}H_\eps)^\ast\Vert$. Indeed, the matrix $H_\eps^{-1}\mathcal{A}H_\eps-(H_\eps^{-1}\mathcal{A}H_\eps)^\ast$ is block diagonal with $m$ blocks of the type $H^{-1}\mathcal{A}H-(H^{-1}\mathcal{A}H)^\ast$. This is the type of matrix which has been estimated in Subsection 4.3 in \cite {GR:11}. In detail, $\Vert H^{-1}\mathcal{A}H-(H^{-1}\mathcal{A}H)^\ast\Vert \le c_3\eps^\alpha\lara{\xi}$ and therefore
\beq
\label{est_3_c1}
\Vert H_\eps^{-1}\mathcal{A}H_\eps-(H_\eps^{-1}\mathcal{A}H_\eps)^\ast\Vert\le c_3\eps^\alpha\lara{\xi}.
\eeq

Finally, if we consider now the matrix of the lower order terms $H_\eps^{-1}\mathcal{L}H_\eps-(H_\eps^{-1}\mathcal{L}H_\eps)^\ast$ we easily sees that it is made of $m$ blocks of the type $({\rm det}H)^{-1}$ times a matrix with $0$-order symbols bounded with respect to $\eps$ (see Subsection 4.4. in \cite{GR:11}). More precisely by following the arguments of Proposition 17(iv) in \cite{GR:11} we get the estimate
\beq
\label{est_4_c1}
\Vert H_\eps^{-1}\mathcal{L}H_\eps-(H_\eps^{-1}\mathcal{L}H_\eps)^\ast\Vert\le c_4\eps^{\alpha(1-m)}.
\eeq

We now insert \eqref{est_1_c1}, \eqref{est_2_c1}, \eqref{est_3_c1} and \eqref{est_4_c1} in the energy estimate \eqref{energy}. We obtain
\begin{multline}
\label{energy_1}
\partial_t |W(t,\xi)|^2\le 2(\rho'(t)\lara{\xi}^{\frac{1}{s}}+c_1\eps^{-1}+c_2\eps^{-1}+c_3\eps^{\alpha}\lara{\xi}+c_4\eps^{\alpha(1-m)})|W(t,\xi)|^2\\
\le (2\rho'(t)\lara{\xi}^{\frac{1}{s}}+C_1\eps^{-1}+C_2\eps^{\alpha}\lara{\xi}+C_3\eps^{\alpha(1-m)})|W(t,\xi)|^2.
\end{multline}
We now set $\eps=\lara{\xi}^{-\gamma}$ in \eqref{energy_1} and we compare the terms
\[
\lara{\xi}^{\gamma},\quad \lara{\xi}^{1-\gamma\alpha},\quad \lara{\xi}^{\gamma\alpha(m-1)}.
\]
For  $\gamma=\min\{\frac{1}{1+\alpha},\frac{1}{\alpha m}\}$ one has that 
\[
\max\{\gamma,\gamma\alpha(m-1)\}\le 1-\gamma\alpha
\]
and therefore
\[
\partial_t |W(t,\xi)|^2\le (2\rho'(t)\lara{\xi}^{\frac{1}{s}}+C\lara{\xi}^{1-\gamma\alpha})|W(t,\xi)|^2.
\]

\subsection{Conclusion of the proof of Theorem \ref{thm1}.}
Let $\rho(t)=\rho(0)-\kappa t$, where $\kappa>0$. If 
\beq
\label{s_c1}
\frac{1}{s}>1-\gamma\alpha=1-\min\biggl\{\frac{\alpha}{1+\alpha},\frac{1}{m}\biggr\}=\max\biggl\{\frac{1}{1+\alpha},\frac{m-1}{m}\biggr\},
\eeq
for $|\xi|$ large enough we have that $\partial_t |W(t,\xi)|^2\le 0$, i.e., $W(t,\xi)=W(0,\xi)$. Therefore,
 
\begin{multline}
\label{last_estimate}
|V(t,\xi)|
=\esp^{-\rho(t)\lara{\xi}^{\frac{1}{s}}}\frac{1}{\det H_\eps(t,\xi)}|H_\eps(t,\xi)||W(t,\xi)|\le \\
\esp^{-\rho(t)\lara{\xi}^{\frac{1}{s}}}\frac{1}{\det H_\eps(t,\xi)}|H_\eps(t,\xi)||W(0,\xi)|=\\
\esp^{(-\rho(t)+\rho(0))\lara{\xi}^{\frac{1}{s}}}\frac{\det H_\eps(0,\xi)}{\det H_\eps(t,\xi)}|H_\eps(t,\xi)||H_\eps^{-1}(0,\xi)||V(0,\xi)|,
\end{multline}
where, arguing on the block level for $\gamma$ as above, we have
\[
\frac{\det H_\eps(0,\xi)}{\det H_\eps(t,\xi)}|H_\eps(t,\xi)||H_\eps^{-1}(0,\xi)|\le c\,\eps^{-\alpha\frac{(m-1)m}{2}}=c\lara{\xi}^{\gamma\alpha\frac{(m-1)m}{2}}.
\]
It follows that,
\[
|V(t,\xi)|\le c\esp^{\kappa T\lara{\xi}^{\frac{1}{s}}}\lara{\xi}^{\gamma\alpha\frac{(m-1)m}{2}}|V(0,\xi)|.
\]
By choosing $\kappa$ small enough we can conclude that $|V(t,\xi)|\le c'\esp^{-\delta\lara{\xi}^{\frac{1}{s}}}$ for some $c',\delta>0$. By the Paley-Wiener characterisation of Gevrey functions this yields to the existence and uniqueness of the solution $U\in C^1([0,T];\gamma^s(\R^n))$ of the Cauchy problem \eqref{CP_syst_Syl} and therefore to the Gevrey well-posedness of the original Cauchy problem \eqref{CP_syst_fin}.


\begin{remark}
\label{rem_Yu}
Note that \eqref{s_c1} implies
\[
s<1+\min\biggl\{\alpha,\frac{1}{m-1}\biggr\}.
\]

This is an improvement in terms of Gevrey order of Yuzawa's  and Kajitani's result in \cite{KY:06, Yu:05}. Indeed, Yuzawa first for t-dependent systems (without lower order terms) in \cite{Yu:05} and later Yuzawa and Kajitani for $(t,x)$-dependent systems in \cite {KY:06} have proven well-posedness in the Gevrey class $\gamma^s$, with 
\[
1\le s<1+\frac{\alpha}{m}.
\]
It is easy to see that
\beq
\label{Yu_GR}
\frac{\alpha}{m}\le \min\biggl\{\alpha,\frac{1}{m-1}\biggr\}.
\eeq
 \end{remark}

\begin{remark}
\label{rem_ultra}
The strategy adopted in the proof of Theorem \ref{thm1} shows how the energy estimate used for scalar equations in \cite{GR:11} can be directly applied to systems after reduction to block Sylvester form to obtain Gevrey well-posedness. In the same way one can get well-posedness in spaces of ultradistributions. In other words, Theorem \ref{thm2} is proven by arguing on the reduced Cauchy problem \eqref{CP_syst_Syl} as in Subsection 4.5 from the aforementioned
paper. 
\end{remark}


\begin{thebibliography}{CDGS79}

\bibitem[CDGS79]{Colombini-deGiordi-Spagnolo-Pisa-1979}
F.~Colombini, E.~De~Giorgi, and S.~Spagnolo.
\newblock Sur les {\'e}quations hyperboliques avec des coefficients qui ne
  d{\'e}pendent que du temps.
\newblock {\em Ann. Scuola Norm. Sup. Pisa Cl. Sci. (4)}, 6(3):511--559, 1979.

\bibitem[CJS87]{Colombini-Jannelli-Spagnolo:Annals-low-reg}
F.~Colombini, E.~Jannelli, and S.~Spagnolo.
\newblock Nonuniqueness in hyperbolic {C}auchy problems.
\newblock {\em Ann. of Math. (2)}, 126(3):495--524, 1987.

\bibitem[CK02]{ColKi:02}
F.~Colombini and T.~Kinoshita.
\newblock On the {G}evrey well posedness of the {C}auchy problem for weakly
  hyperbolic equations of higher order.
\newblock {\em J. Differential Equations}, 186(2):394--419, 2002.

\bibitem[CS82]{Colombini-Spagnolo:Acta-ex-weakly-hyp}
F.~Colombini and S.~Spagnolo.
\newblock An example of a weakly hyperbolic {C}auchy problem not well posed in
  {$C^{\infty }$}.
\newblock {\em Acta Math.}, 148:243--253, 1982.

\bibitem[DS98]{DS}
P.~D'Ancona and S.~Spagnolo.
\newblock Quasi-symmetrization of hyperbolic systems and propagation of the
  analytic regularity.
\newblock {\em Boll. Unione Mat. Ital. Sez. B Artic. Ric. Mat. (8)},
  1(1):169--185, 1998.

\bibitem[G15]{G:15}
C.~Garetto.
\newblock On hyperbolic equations and systems with non-regular time dependent coefficients
\newblock {\em J. Differential Equations}, http://arxiv.org/abs/1504.03716, 2015.

\bibitem[GR12]{GR:11}
C.~Garetto and M.~Ruzhansky.
\newblock On the well-posedness of weakly hyperbolic equations with
  time-dependent coefficients.
\newblock {\em J. Differential Equations}, 253(5):1317--1340, 2012.

\bibitem[GR13]{GR:12}
C.~Garetto and M.~Ruzhansky.
\newblock Weakly hyperbolic equations with non-analytic coefficients and lower
  order terms.
\newblock {\em Math. Ann.}, 357(2):401--440, 2013.

\bibitem[GR14]{Garetto-Ruzhansky:JMAA-2014}
C.~Garetto and M.~Ruzhansky.
\newblock A note on weakly hyperbolic equations with analytic principal part.
\newblock {\em J. Math. Anal. Appl.}, 412(1):1--14, 2014.

\bibitem[GR15]{Garetto-Ruzhansky:ARMA}
C.~Garetto and M.~Ruzhansky.
\newblock Hyperbolic second order equations with non-regular time dependent
  coefficients.
\newblock {\em Arch. Ration. Mech. Anal.}, 217(1):113--154, 2015.

\bibitem[Jan84]{Jannelli:analytic-CPDE-1984}
E.~Jannelli.
\newblock Linear {K}ovalevskian systems with time-dependent coefficients.
\newblock {\em Comm. Partial Differential Equations}, 9(14):1373--1406, 1984.

\bibitem[Kaj86]{Kajitani:analytic-CPDE-1986}
K.~Kajitani.
\newblock Global real analytic solutions of the {C}auchy problem for linear
  partial differential equations.
\newblock {\em Comm. Partial Differential Equations}, 11(13):1489--1513, 1986.

\bibitem[KS06]{KS}
T.~Kinoshita and S.~Spagnolo.
\newblock Hyperbolic equations with non-analytic coefficients.
\newblock {\em Math. Ann.}, 336(3):551--569, 2006.

\bibitem[KY06]{KY:06}
K.~Kajitani and Y.~Yuzawa.
\newblock The {C}auchy problem for hyperbolic systems with {H}{\"o}lder
  continuous coefficients with respect to the time variable.
\newblock {\em Ann. Sc. Norm. Super. Pisa Cl. Sci. (5)}, 5(4):465--482, 2006.

\bibitem[Yuz05]{Yu:05}
Y.~Yuzawa.
\newblock The {C}auchy problem for hyperbolic systems with {H}{\"o}lder
  continuous coefficients with respect to time.
\newblock {\em J. Differential Equations}, 219(2):363--374, 2005.

\end{thebibliography}

\end{document}